\documentclass{amsart}

\usepackage{amssymb}
\usepackage[all]{xy}
\CompileMatrices

\newtheorem*{introthm}{Theorem}
\newtheorem{theorem}{Theorem}[section]
\newtheorem{corollary}[theorem]{Corollary}
\newtheorem{proposition}[theorem]{Proposition}

\newtheorem{lemma}[theorem]{Lemma}

\theoremstyle{definition}
\newtheorem{definition}[theorem]{Definition}
\newtheorem{remark}[theorem]{Remark}

\def\quot{/\!\!/}
\def\stern{\! * \!}
\def\mal{\! \cdot \!}
\def\rq#1{\widehat{#1}}
\def\t#1{\widetilde{#1}}
\def\b#1{\overline{#1}}
\def\KK{{\mathbb K}}
\def\CC{{\mathbb C}}

\def\ZZ{{\mathbb Z}}
\def\NN{{\mathbb N}}
\def\QQ{{\mathbb Q}}
\def\PP{{\mathbb P}}
\def\SL{{\rm SL}}

\def\Div{{\rm div}}
\def\WDiv{{\rm WDiv}}
\def\CDiv{{\rm CDiv}}
\def\Supp{{\rm Supp}}
\def\Char{{\rm Char}}
\def\Spec{{\rm Spec}}

\def\pr{{\rm pr}}
\def\id{{\rm id}}

\begin{document}

\title[GIT based on Weil divisors]
      {Geometric Invariant Theory \\
       based on Weil divisors}  
\author{J\"urgen Hausen} 
\email{hausen@mfo.de}
\address{Mathematisches Forschungsinstitut Oberwolfach \\
         Schwarzwaldstra\ss e 9-11 \\ 
         77709 Oberwolfach-Wal\-ke \\ 
         Germany}

\begin{abstract}
Given an action of a reductive group on a normal 
variety, we describe all invariant open subsets 
admitting a good quotient with a quasiprojective 
or a divisorial quotient space.
We obtain several new Hilbert-Mumford type theorems,
and we extend a projectivity criterion 
of Bia\l ynicki-Biru\-la and \'Swi\c{e}cicka for 
varieties with semisimple group action 
from the smooth to the singular case.
\end{abstract}

\maketitle

\section*{Introduction}

This article is devoted to a central task
of Geometric Invariant Theory, formulated 
in~\cite{BB}: 
Given an action of a reductive group $G$ on a normal
variety $X$, describe all $G$-invariant open 
subsets $U \subset X$ admitting a good quotient,
that means a $G$-invariant affine morphism
$U \to U \quot G$ such that the structure sheaf
of $U \quot G$ equals the sheaf of invariants 
$p_{*}(\mathcal{O}_{U})^{G}$.
We call these $U$ for the moment the 
$G$-sets.

In~\cite{Mu}, Mumford obtains $G$-sets 
with quasiprojective quotient spaces.
Given a $G$-linearized line bundle
$L \to X$, that means that $G$ acts
on the total space making the projection
equivariant and inducing linear maps
on the fibres,
he calls a point $x \in X$ semistable,
if some positive power of $L$ admits a 
$G$-invariant section $f$ such that removing 
the zeroes gives an affine neighbourhood $X_f$ 
of~$x$.

The set $X^{ss}(L)$ of semistable 
points of a $G$-linearized line bundle $L$
admits a good quotient 
$X^{ss}(L) \to X^{ss}(L) \quot G$ 
with a quasiprojective quotient space.
For smooth $X$, basically 
all quasiprojective quotient spaces
arise in this way:
every $G$-set $U$ with $U \quot G$
quasiprojective is $G$-saturated in some 
$X^{ss}(L)$, that means that $U$ is saturated 
with respect to the quotient map.

For singular $X$, Mumford's 
method does not provide all 
quasiprojective quotients,
see Example~\ref{counterexample}.
Here, replacing the bundles $L$ with 
Weil divisors $D$ yields a more 
rounded picture: we define a
$G$-linearization of $D$ to be 
a certain lifting of the $G$-action to
$\Spec(\mathcal{A})$, 
where $\mathcal{A} = \oplus_{n \ge 0} \mathcal{O}(nD)$,
and the set $X^{ss}(D)$ of semistable points
is the union of all affine sets
$X_f$, where $f \in \mathcal{A}(X)$ is $G$-invariant 
and homogeneous of positive degree.  
The first result is Theorem~\ref{qpquots}:

\begin{introthm}
Let a reductive group $G$ act on a normal variety $X$.
\begin{enumerate}
\item For any $G$-linearized Weil divisor $D$ on $X$, there is
  a good quotient $X^{ss}(D) \to X^{ss}(D) \quot G$ with 
  a quasiprojective variety $X^{ss}(D) \quot G$.
\item If $U \subset X$ is open, $G$-invariant, and has a good 
  quotient $U \to U \quot G$ with $U \quot G$ quasiprojective,
  then $U$ is a $G$-saturated subset of some $X^{ss}(D)$.
\end{enumerate}
\end{introthm}

However, the quasiprojective quotient spaces 
are by far not the whole story, and a further 
aim is to complement also the picture developed
in~\cite{Ha0}
for {\em divisorial\/} quotient spaces, 
that means (possibly nonseparated)
prevarieties $Y$ such that
every $y \in Y$ has an affine neighbourhood 
$Y \setminus \Supp(E)$ with an
effective Cartier divisor~$E$,
see~\cite{Bo} and~\cite{SGA}.
For the occurence of nonseparatedness  
in quotient constructions, compare
also~\cite[Prop.~1.9]{Mu}, \cite[Cor.~1.3]{Su},
and~\cite{AcHa0}.

To obtain divisorial quotient spaces, we work 
with finitely generated groups $\Lambda$ 
of Weil divisors.
Similarly as before, $G$-linearization of 
such a $\Lambda$ is a lifting of the $G$-action
to $\Spec(\mathcal{A})$, where now
$\mathcal{A} = \oplus_{D \in \Lambda} 
\mathcal{O}(D)$. 
We also have a notion of semistability, 
and the resulting statements generalize~\cite{Ha0}, 
see~Theorem~\ref{divisorialquotients}:

\begin{introthm}
Let a reductive group $G$ act on a normal variety $X$.
\begin{enumerate}
\item For any $G$-linearized group $\Lambda$ of Weil divisors 
  on $X$, there is a good quotient 
  $X^{ss}(\Lambda) \to X^{ss}(\Lambda) \quot G$ with 
  a divisorial prevariety $X^{ss}(\Lambda) \quot G$.
\item If $U \subset X$ is open, $G$-invariant, and has a good 
  quotient $U \to U \quot G$ with $U \quot G$ divisorial,
  then $U$ is a $G$-saturated subset of some $X^{ss}(\Lambda)$.
\end{enumerate}
\end{introthm}

A simple example shows that, in general, 
the respective sets of semistable points of 
a single linearized divisor $D$ and the group 
$\ZZ D$ differ, compare~\cite[Example~3.5]{Ha0}.
Let $G := \CC^*$ act linearly on $X := \CC^2$ via
$$
t \mal (z,w) := (tz,t^{-1}w).
$$
Consider the invariant divisor $D := {\rm div}(z)$ 
on $X$.
Then $D$ as  well as the group $\Lambda := \ZZ D$
are canonically $G$-linearized, via the induced action
of $G$ on the function field.
According to the respective definitions~\ref{semistabdef}
and~\ref{semistabgroupdef} of semistability,
one obtains
$$
X^{ss}(D) = \CC^* \times \CC,
\qquad
X^{ss}(\Lambda) = \CC^2 \setminus \{(0,0)\}.
$$
Moreover, let us have a look at the quotient spaces.
For the first set, the quotient space is 
the affine line, 
whereas in the second case a true
(divisorial) {\em pre}-variety occurs: 
the affine line with a doubled point.

For practical purposes, it is
often helpful to perform
the construction of $G$-sets
by means of subtori of $G$.
Classically, this is done by the
Hilbert-Mumford Lemma~\cite[Thm.~2.1]{Mu}:
for a $G$-linearized ample bundle $L$
on a projective variety $X$, 
it gives a semistability criterion in 
terms of one parameter subgroups;
here, we deal with the following version,
involving a maximal torus 
$T  \subset G$, compare~\cite{BB} 
and~\cite{AS}:
$$
X^{ss}(L,G) = \bigcap_{g \in G} g \mal X^{ss}(L,T).
$$

In this form, the statement allows 
a far reaching generalization; 
in particular, the hypotheses
of projectivity and ampleness can be dropped,
see Theorem~\ref{himusemistab}:

\begin{introthm}
Let a reductive group $G$ act on a normal 
variety $X$, and let $T \subset G$ be a maximal torus.
\begin{enumerate} 
\item Let $D$ be a $G$-linearized Weil 
divisor on $X$. Then we have:
$$ X^{ss}(D,G) = \bigcap_{g \in G} g \mal X^{ss}(D,T). $$
\item Let 
$\Lambda \subset \WDiv(X)$ be a $G$-linearized 
subgroup. Then we have:
$$ 
X^{ss}(\Lambda,G) 
= 
\bigcap_{g \in G} g \mal X^{ss}(\Lambda,T). 
$$
\end{enumerate}
\end{introthm}

Finally, in Section~\ref{section7},
we focus on the case of a semisimple group $G$.
We ask for maximal $G$-sets,
compare~\cite{BB}:
A {\em qp-maximal $G$-set\/} is a
$G$-set $U \subset X$ with $U \quot G$ 
quasiprojective such that $U$ does 
not occur as a $G$-saturated 
proper subset in some 
$U' \subset X$ with the same properties. 
Similarly, a {\em d-maximal $G$-set\/} is 
a subset having the analogous properties 
with respect to divisorial
quotient spaces.

Reducing the construction of these 
sets to the construction of the qp- and 
the d-maximal $T$-sets for 
a maximal torus $T \subset G$
amounts to tackling Bia\l ynicki-Birula's 
Conjecture~\cite[12.1]{BB}:
Given a maximal $T$-set $U \subset X$ which is invariant
under the normalizer $N \subset G$ of $T$, he asks 
if the following set is open and admits a good quotient 
by $G$:
$$ W(U) := \bigcap_{g \in G} g \mal U .$$

Here are the known positive results 
concerning qp- and d-maximal $T$-sets 
$U \subset X$:
The case of $G = \SL_{2}$ acting on a smooth $X$
is settled in~\cite[Thm.~9]{BBSw2} 
and~\cite[Thm.~2.2]{Ha2}.
If $U \quot T$ is projective and $X$ is smooth,
then~\cite[Cor.~1]{BBSw3} gives positive
answer for a general connected semisimple 
group $G$. 
Moreover, the problem is solved in the case $U = X$, 
see~\cite[Thm.~12.4]{BB} and~\cite[Thm.~5.1]{Ha0}.
We show in Theorem~\ref{semisimhimu} and 
Corollary~\ref{semisimple2hilbmum}:

\begin{introthm}
Let $G$ be a connected semisimple group, and 
$T \subset G$ a maximal torus with normalizer 
$N \subset G$.
Let $X$ be a normal $G$-variety, $U \subset X$ 
an $N$-invariant open subset, 
and $W(U)$ the intersection of all
translates $g \mal U$, where $g \in G$.
\begin{enumerate}
\item If $U \subset X$ is a qp-maximal $T$-set,
   then $W(U)$ is open and $T$-saturated in $U$,
   and there is a good quotient $W(U) \to W(U) \quot G$
   with $W(U) \quot G$ quasiprojective.
\item If $U$ admits a good quotient $U \to U \quot T$
   with $U \quot T$ projective, then $W(U)$ is open 
   and $T$-saturated in $U$, and there is a good 
   quotient $W(U) \to W(U) \quot G$ with $W(U) \quot G$
   projective.
\item If $U \subset X$ is a d-maximal $N$-set,
   then $W(U)$ is open and $T$-saturated in $U$,
   and there is a good quotient $W(U) \to W(U) \quot G$
   with  $W(U) \quot G$ divisorial.
\end{enumerate}
\end{introthm}

In the setting of~(ii), we can prove much more.
It turns out that $U$ and $W(U)$ are the sets of 
semistable points of an ordinary linearized ample line bundle, 
and  --- even more surprising --- that $X$ is projective.
This extends the main result of~\cite{BBSw3}
from the smooth to the normal case and thus gives 
an answer to the problem discussed 
in~\cite[Remark p.~965]{BBSw3}.
More precisely, we prove in Theorem~\ref{projcrit}:

\begin{introthm}
Let $G$ be a connected semisimple group,
$T \subset G$ a maximal torus with normalizer 
$N \subset G$, and $X$ be a normal $G$-variety.
Suppose that $U \subset X$ is $N$-invariant, open 
and admits a good quotient 
$U \to U \quot T$ with $U \quot T$ projective.
Then there is an ample $G$-linearized line bundle 
$L$ on $X$ with $U = X^{ss}(L,T)$,
we have $X = G \mal U$, and $X$ is projective.
\end{introthm}

I would like to thank the referee for his 
helpful suggestions and comments on the 
earlier versions of this article.

\section{Polyhedral semigroups and $G$-linearization}
\label{section2}

In this section, we transfer Mumford's concepts
of~\cite[Sec.~1.3]{Mu} to the framework of
Weil divisors.
We introduce polyhedral semigroups of
Weil divisors, and define the notion of 
a $G$-linearization for such a semigroup.
Moreover, we give a geometric interpretation
of this concept, and provide basic statements 
concerning existence and uniqueness of 
linearizations.

Throughout the whole article, we work over an 
algebraically closed field $\KK$ of 
characteristic zero.
In this section, $X$ denotes an irreducible 
normal prevariety over $\KK$,
that means that $X$ is an integral, normal,
but possibly nonseparated scheme of finite 
type over $\KK$, 
compare also~\cite[Sec.~I.2.2]{Hu}.
The word ``point'' always refers to a closed point.
 
By $\WDiv(X)$ we denote the group of Weil 
divisors of $X$,
and $\CDiv(X) \subset \WDiv(X)$ is the 
subgroup of Cartier divisors.
For a finitely generated subsemigroup 
$\Lambda \subset \WDiv(X)$, 
let $\Gamma(\Lambda) \subset \WDiv(X)$ denote 
the subgroup generated by $\Lambda$.
We say that the semigroup $\Lambda$ 
is {\em polyhedral}, if 
it is the intersection of $\Gamma(\Lambda)$ 
with a convex polyhedral cone
in $\QQ \otimes_{\ZZ} \Gamma(\Lambda)$.

Fix a polyhedral semigroup
$\Lambda \subset \WDiv(X)$.
Since we assumed $X$ to be normal, 
there is an associated $\mathcal{O}_{X}$-module
$\mathcal{O}_{X}(D)$ of rational functions
for any $D \in \Lambda$. 
In fact, multiplication in the function field 
$\KK(X)$ gives even rise to a $\Lambda$-graded 
$\mathcal{O}_{X}$-algebra:
$$ 
\mathcal{A} 
:= \bigoplus_{D \in \Lambda} \mathcal{A}_{D}
:= \bigoplus_{D \in \Lambda} \mathcal{O}_{X}(D). 
$$

Now, let $G$ be a linear algebraic group,
and let $G$ act on $X$.
That means in particular that this action 
is given by a morphism
$\alpha \colon G \times X \to X$,
and, denoting by $\mu \colon G \times G \to G$
the multiplication map, we have 
a commutative diagram
$$
\xymatrix{
G \times G \times X 
\ar[rr]^{\id_G \times \alpha} 
\ar[d]_{\mu \times \id_X}
& &
G \times X 
\ar[d]^{\alpha}
\\
G \times X \ar[rr]_{\alpha} 
& &
X
}
$$
Similarly to~\cite[Def.~1.6]{Mu},
the definition of a $G$-linearization of 
the semigroup $\Lambda \subset \WDiv(X)$
is formulated in terms of $\mathcal{A}$,
the above maps and the projection maps
\begin{eqnarray*}
\pr_{G \times X} \colon G \times G \times X \to G \times X,
& \qquad &
(g_1,g_2,x) \mapsto (g_2,x),
\\
\pr_X \colon G \times X \to X, 
& \qquad & (g,x) \mapsto x.
\end{eqnarray*}

\begin{definition}
\label{glindef}
A {\em $G$-linearization\/} of $\Lambda$ is an isomorphism 
$\Phi \colon \alpha^* \mathcal{A} \to \pr_X^* \mathcal{A}$ 
of $\Lambda$-graded $\mathcal{O}_{G \times X}$-algebras
such that $\Phi$ is the identity in degree zero, 
and the following diagram is commutative:
$$
\xymatrix{
(\id_G \times \alpha)^* \alpha^* \mathcal{A}
\ar[rr]^{(\id_G \times \alpha)^* \Phi}
& &
(\id_G \times \alpha)^* \pr_X^* \mathcal{A}
\ar@{=}[rr]
& &
{\pr_{G \times X}^* \alpha^* \mathcal{A}}
\ar[d]^{\pr_{G \times X}^* \Phi}
\\
(\mu \times \id_X)^* \alpha^* \mathcal{A}
\ar[rr]_{(\mu \times \id_X)^* \Phi}
\ar@{=}[u]
& &
(\mu \times \id_X)^* \pr_X^* \mathcal{A}
\ar@{=}[rr]
& &
{\pr_{G \times X}^* \pr_X^* \mathcal{A}}
}
$$
\end{definition}

Note that if $\Lambda = \oplus_{n \ge 0} nD$ 
with a Cartier divisor $D$, then 
the $G$-linearizations
$\Phi \colon \alpha^* \mathcal{A} \to \pr_X^* \mathcal{A}$ 
of $\Lambda$ correspond to the $G$-linearizations
of the invertible sheaf $\mathcal{O}_X(D)$ 
in the sense of~\cite[Def.~1.6]{Mu} via
passing to the corresponding map in degree one
$\Phi_1 \colon \alpha^* \mathcal{O}_X(D) \to \pr_X^* \mathcal{O}_X(D)$.

In order to interprete Def.~\ref{glindef} geometrically,
look at the scheme 
$\t{X} := \Spec(\mathcal{A})$
over $X$.
Note that the $\Lambda$-grading of $\mathcal{A}$ 
defines an action of the torus 
$S := \Spec(\KK[\Gamma(\Lambda)])$ on $\t{X}$.
We list some properties; for example,
over the smooth locus,
the canonical map $q \colon \t{X} \to X$ 
is locally trivial with an
affine toric variety as fibre:

\begin{proposition}\label{toricbundle}
Let $U \subset X$ be an open subset such that
every $D \in \Lambda$ is Cartier on $U$, and set 
$\t{U} := q^{-1}(U)$. 
\begin{enumerate}
\item The map $q \colon  \t{U} \to U$ is locally trivial with 
  typical fibre $\t{U}_{x} \cong \Spec(\KK[\Lambda])$.
  The open set $\rq{U} \subset \t{U}$ of free $S$-orbits is an
  $S$-principal bundle over $U$.
\item The inclusion $\rq{U} \subset \t{U}$ corresponds to the 
  inclusion $\mathcal{A} \subset \mathcal{B}$ of the graded
  $\mathcal{O}_{U}$-algebras $\mathcal{A}$ and $\mathcal{B}$
  arising from $\Lambda$ and $\Gamma(\Lambda)$.
\item For any homogeneous section $f \in \mathcal{A}(U)$, its zero set 
  as a function on $\rq{U}$ equals the set 
  $\rq{U} \cap q^{-1}(\Supp(\Div(f)+D))$.
\end{enumerate}
\end{proposition}

\begin{proof}
Consider the group $\Gamma(\Lambda)$ generated 
by $\Lambda$ and its $\mathcal{O}_{X}$-algebra~$\mathcal{B}$.
Locally, $\mathcal{B}$ is a Laurent monomial algebra 
over $\mathcal{O}_{U}$, i.e., for small affine open $V \subset U$, 
we have a graded isomorphism over $\mathcal{O}(V)$:
$$ 
\mathcal{B}(V) 
\; \cong \; 
\mathcal{O}(V) \otimes_{\KK} \KK[\Gamma(\Lambda)].
$$
Cutting down this to the subsemigroup $\Lambda \subset \Gamma(\Lambda)$
and the associated subalgebra $\mathcal{A} \subset \mathcal{B}$, we obtain
local triviality of $q \colon \t{U} \to U$. The remaining statements follow 
then easily.
\end{proof}

Any $G$-linearization 
$\Phi \colon \alpha^* \mathcal{A} \to \pr_X^* \mathcal{A}$
of the polyhedral semigroup 
$\Lambda \subset \WDiv(X)$
defines a commutative diagram
\begin{equation}
\label{spec}
\vcenter{
\xymatrix{
{\Spec(\pr_X^*\mathcal{A})}
\ar[rr]^{\Spec(\Phi)}
\ar@{=}[d]
& &
{\Spec(\alpha^* \mathcal{A})}
\ar[rr]
& &
{\Spec(\mathcal{A})}
\ar@{=}[d]
\\
G \times \t{X}
\ar[rrrr]_{\t{\alpha}}
& & & &
{\t{X}}
}}
\end{equation}
Note that $\Spec(\alpha^{*} \mathcal{A})$
is the fibre product of $\alpha \colon G \times X \to X$ 
and the canonical map $\t{X} \to X$.
Then the right upper arrow is merely the 
projection to $\t{X}$. 

\begin{proposition}
\label{liftedaction}
\begin{enumerate}
\item The map 
$\t{\alpha} \colon G \times \t{X} \to \t{X}$ 
is a $G$-action that  
commutes with the $S$-action on $\t{X}$,
and makes the canonical map 
$\t{X} \to X$ equivariant. 
\item For every action 
$\t{\alpha} \colon G \times \t{X} \to \t{X}$ 
as in~(i), there is a unique $G$-linearization 
$\Phi \colon \alpha^* \mathcal{A} \to \pr_X^* \mathcal{A}$
making the diagram~(\ref{spec}) commutative.
\end{enumerate}
\end{proposition}

\begin{proof}
For (i), note that $q \circ \t{\alpha}$ 
equals $\alpha \circ (\id_G \times q)$, 
because $\Phi$ is the identity in degree 
zero. 
Moreover, the commutative diagram of Def.~\ref{glindef} 
yields the associativity law 
of a group action for $\t{\alpha}$,
and $e_G  \in G$ acts trivially
because $\Phi$ is an isomorphism.
Finally, the actions of $G$ and $S$ commute, 
because $\t{\alpha}$ has graded
comorphisms.

To verify~(ii), we use that $\Spec(\alpha^{*} \mathcal{A})$
is the fibre product of $\alpha \colon G \times X \to X$ 
and $\t{X} \to X$.
By the universal property, 
$\t{\alpha} \colon G \times \t{X} \to \t{X}$ 
lifts to a unique morphism 
$\Spec(\pr_{X}^{*} \mathcal{A})
\to \Spec(\alpha^{*} \mathcal{A})$.
It is straightforward to check that 
this morphism stems from a $G$-linearization
$\Phi \colon \alpha^* \mathcal{A} \to \pr_X^* \mathcal{A}$.
\end{proof}

Eventually, via the lifted $G$-action on $\t{X}$,
we associate to any $G$-linearization of $\Lambda$
a {\em graded $G$-sheaf structure\/} on 
$\mathcal{A}$.
The latter is a collection of graded 
$\mathcal{O}(U)$-algebra homomorphisms 
$\mathcal{A}(U) \to \mathcal{A}(g \mal U)$,
$f \mapsto g \mal f$,
being compatible with group operations in $G$ 
and with restriction and algebra operations
in $\mathcal{A}$;
thereby $G$ acts as usual on the structure sheaf 
$\mathcal{O}_{X}$ via $g \mal f(x) := f(g^{-1} \mal x)$.
 
\begin{proposition}
\label{gsheaf}
Let $\Phi \colon \alpha^* \mathcal{A} \to \pr_X^* \mathcal{A}$
be a $G$-linearization.
Then there is a unique graded $G$-sheaf structure on 
$\mathcal{A}$ satisfying
$g \mal f(\t{x}) := f(g^{-1} \mal \t{x})$
for any $\t{x} \in \t{X }$ lying 
over the smooth locus of $X$.
For every $G$-invariant open $U \subset X$, the
induced representation of $G$ on $\mathcal{A}(U)$ is 
rational.
\end{proposition} 

\begin{proof}
Over the smooth locus of $X$, 
we may define the $G$-sheaf structure
according to
$g \mal f(\t{x}) := f(g^{-1} \mal \t{x})$.
By normality, it uniquely extends to $X$.
Rationality of the induced representations
follows, e.g., from~\cite[Lemma~2.5]{DMV}.
\end{proof}

\begin{remark}
Let $\Phi \colon \alpha^* \mathcal{A} \to \pr_X^* \mathcal{A}$
be a $G$-linearization.
Then a section $f \in \mathcal{A}(X)$ is 
invariant with respect to the induced
$G$-representation on $\mathcal{A}(X)$ if 
and only if $\Phi(\alpha^*(f)) = \pr_X^*(f)$
holds.
\end{remark}

We give two existence statements for $G$-linearizations.
The first one is the analogue of Mumford's 
result~\cite[Cor.~1.6]{Mu} and~\cite[Prop.~2.4]{DMV}.
We use the following terminology:
Given polyhedral semigroups
$\Lambda' \subset \Lambda$,
we say that $\Lambda'$ {\em is of finite index\/} 
in $\Lambda$ if there is a positive $n \in \ZZ$ with 
$n \Lambda \subset \Lambda'$.

\begin{proposition}\label{existence1}
Suppose that $X$ is separated and that $G$ is connected.
Then, for any polyhedral semigroup $\Lambda \subset \WDiv(X)$,
some subsemigroup $\Lambda' \subset \Lambda$
of finite index admits a $G$-linearization.
\end{proposition}

\begin{proof} 
By normality, it suffices to provide 
a $G$-linearization of $\mathcal{A}$ over
the smooth locus. Hence, we may assume that
$\Gamma \subset \CDiv(X)$ holds.
Consider the group 
$\Gamma(\Lambda) \subset \CDiv(X)$ 
generated by $\Lambda$, and fix any basis
$D_1, \ldots, D_k$ of $\Gamma(\Lambda)$.
Then~\cite[Prop.~2.4]{DMV} gives us 
$n_i \ge 1$ and linearizations in the sense
of~\cite[Def.~1.6]{Mu}:
$$
\Phi_i \colon 
\alpha^* \mathcal{O}_X(n_i D_i) 
\to 
\pr_X^* \mathcal{O}_X(n_i D_i).
$$

Let $\Gamma' \subset \Gamma(\Lambda)$
be the subgroup generated by the $n_iD_i$.
Then, via tensoring the $\Phi_i$, we obtain 
for each $D \in \Gamma' \cap \Lambda$ 
an isomorphism
$\alpha^* \mathcal{A}_D \to \pr_X^* \mathcal{A}_D$.
These maps are compatible with the multiplicative
structures of $\alpha^*\mathcal{A}$ and 
$\pr_X^*\mathcal{A}$, and hence fit together to 
a linearization of $\Gamma' \cap \Lambda$.
\end{proof}

The second existence statement provides
{\em canonical $G$-linearizations\/}.
As usual, we say that a Weil divisor 
$D = \sum n_{E}E$ is {\em $G$-invariant\/} 
if $n_{g \mal E} = n_{E}$ holds for any prime divisor $E$.
The support of a $G$-invariant Weil divisor is
$G$-invariant, whereas its components may be permuted.

Moreover, we have to consider
pullbacks of $\mathcal{A}$
under dominant maps $p \colon Z \to X$,
where $Z$ is normal.
If the inverse image 
$p^{-1}(X')$ of the smooth locus 
$X' \subset X$ has a complement
of codimension at least two in $Z$,
then the pullback
$\CDiv(X') \to \CDiv(p^{-1}(X'))$
induces a map 
$p^* \colon \WDiv(X) \to \WDiv(Z)$,
and we obtain
$$
p^* \mathcal{A}
\; = \;
\bigoplus_{D \in \Lambda} \mathcal{O}_{Z} (p^* D).
$$

\begin{proposition}
\label{existence2}
Let $\Lambda$ consist of $G$-invariant divisors.
Then there is a canonical $G$-linearization
$$ 
\alpha^* \mathcal{A}
\;=\;
\bigoplus_{D \in \Lambda} \mathcal{O}_{G \times X}(\alpha^*D)
\;=\;
\bigoplus_{D \in \Lambda} \mathcal{O}_{G \times X}(\pr_X^*D)
\;=\;
\pr_X^* \mathcal{A}.
$$
The induced $G$-sheaf structure on $\mathcal{A}$ is given
by the usual action of $G$ on the function field
$\KK(X)$ via $g \mal f(x) = f(g^{-1} \mal x)$.
\end{proposition}

\begin{proof} 
We have to show that any $G$-invariant
Weil divisor $D = \sum n_E E$ satisfies $\alpha^*D = \pr_X^* D$.
For this, we consider the isomorphism
$$
\beta \colon G \times X \to G \times X,
\qquad
(g,x) \mapsto (g,g^{-1} \mal x).
$$
Then we have $\beta^* \pr_X^* D = \pr_X^* D$,
and $\pr_X^*D = \beta^* \alpha^* D$.
Since $\beta^*$ has an inverse, the assertion 
follows.   
\end{proof}

We turn to uniqueness properties of $G$-linearizations.
Let $\Char(G)$ denote the group of characters of $G$, i.e., the 
group of all homomorphisms $G \to \KK^{*}$.
For groups $G$ with few characters, we have the following
two statements, compare~\cite[Prop.~1.4]{Mu} and~\cite[Prop.~1.5]{Ha0}:

\begin{proposition}\label{uniqueness}
Let $X$ be separated, and let $\Lambda \subset \WDiv(X)$ 
be a polyhedral semigroup.
\begin{enumerate}
\item If $\Char(G)$ is trivial and $G$ is connected, then any two
  $G$-linearizations of $\Lambda$ coincide.
\item If $\Char(G)$ is finite and $\mathcal{O}^{*}(X) = \KK^{*}$
  holds, then any two $G$-linearizations of $\Lambda$ induce 
  the same $G$-linearization on some $\Lambda' \subset \Lambda$
  of finite index.
\end{enumerate}
\end{proposition}

\begin{proof} 
Again by normality, it suffices to treat
the problem over the smooth locus. 
Then $q \colon \t{X} \to X$ 
is locally trivial with toric fibres, having
$S = \Spec(\KK[\Gamma(\Lambda)])$ as their
big torus.
Given two $G$-linearizations of $\Lambda$,
we denote the two corresponding
$G$-actions on $\t{X}$ by 
$g \mal z$ and $g \stern z$.
Consider the morphism
$$ 
\psi \colon G \times \t{X} \to \t{X},
\qquad
z \mapsto g^{-1} \stern g \mal z. 
$$

For fixed $g$, the map $z \mapsto \psi(g,z)$ is an 
$S$-equivariant bundle automorphism.
Hence, on each fibre it is multiplication 
with an element of the torus $S$.
Consequently, there is a morphism 
$\eta \colon G \times X \to S$ 
such that $\psi$ is of the form
$$ \psi(g,z) = \eta(g,q(z)) \mal z $$ 

In the setting of~(i), Rosenlicht's Lemma~\cite[Lemma~2.1]{FoIv}
yields a decomposition $\eta(g,z) = \chi(g) \beta(q(z))$
with a regular homomorphism $\chi \colon G \to S$ 
and a morphism $\beta \colon X \to S$. 
Since we assumed $G$ to have only trivial characters, 
we can conclude that $\psi$ is the identity map.

If we are in the situation of~(ii), then 
$\mathcal{O}^{*}(X) = \KK^{*}$ implies that
$\psi(g,z) = \chi(g) \mal z$ holds with a 
regular homomorphism $\chi \colon G \to S$.
Hence, after dividing $\t{X}$ by the finite subgroup 
$\chi(G) \subset S$, the two induced $G$-actions coincide.
But this process means replacing $\Lambda$ with a 
subsemigroup of finite index.
\end{proof}

Let us remark that there are simple examples 
showing that for non connected $G$,
one cannot omit the assumption $\mathcal{O}^{*}(X) = \KK^{*}$
in the second statement.

\section{The ample locus}\label{section3}

We introduce the Cartier locus and the 
ample locus of a polyhedral semigroup 
of Weil divisors, and study its behaviour in the case of 
$G$-linearized semigroups.
The considerations of this section prepare the proofs of the 
various Hilbert-Mumford type theorems given later.

Unless otherwise stated, $X$ denotes in this section an 
irreducible normal prevariety. 
Given a polyhedral semigroup $\Lambda \subset \WDiv(X)$,
let $\mathcal{A}$ denote the associated
$\Lambda$-graded $\mathcal{O}_{X}$-algebra.
For a homogeneous local section $f \in \mathcal{A}_{D}(U)$, 
we define its {\em zero set\/} to be
$$ Z(f) := \Supp(\Div(f) + D\vert_{U}). $$

\begin{definition}\label{ampledef}
Let $\Lambda \subset \WDiv(X)$ be a polyhedral 
semigroup with associated $\Lambda$-graded
$\mathcal{O}_{X}$-algebra $\mathcal{A}$.
\begin{enumerate}
\item The {\em Cartier locus\/} of $\Lambda$ is the set 
   of all points $x \in X$ such that every $D \in \Lambda$ 
   is Cartier near $x$.
\item The {\em ample locus\/} of $\Lambda$ is the set 
   of all $x \in X$ admitting an affine neighbourhood 
   $X \setminus Z(f)$ with a homogeneous section 
   $f \in \mathcal{A}(X)$ such that $X \setminus Z(f)$
   is contained in the Cartier locus of $\Lambda$.
\end{enumerate}
\end{definition}

We shall speak of an {\em ample\/} semigroup 
$\Lambda \subset \WDiv(X)$ if the 
ample locus of $\Lambda$ equals $X$.
Thus, ample semigroups consist by definition
of Cartier divisors.
The relations to the usual ampleness concepts 
~\cite{EGA}, \cite{Bo} and~\cite{SGA} are the 
following; recall  
that $X$ is said to be {\em divisorial\/}, if 
every $x \in X$ has an affine neighbourhood 
$X \setminus \Supp(E)$ with an effective
Cartier divisor $E$ on $X$.

\begin{remark}
\begin{enumerate}
\item A polyhedral semigroup of the form 
   $\Lambda = \NN D$ is ample if and only if 
   $D$ is an ample Cartier divisor in the usual sense.
\item An irreducible normal prevariety is divisorial 
  if and only if it admits an ample group of Cartier divisors.
\end{enumerate}
\end{remark}

Let us explain the geometric meaning of the ample locus
of a polyhedral semigroup $\Lambda \subset \CDiv(X)$
in terms of the corresponding toric bundle $q \colon \t{X} \to X$.
Recall from Section~\ref{section2} that $\t{X}$ comes along with 
an action of the torus $S = \Spec(\KK[\Gamma(\Lambda)])$, 
and that the set $\rq{X} \subset \t{X}$ of free $S$-orbits
is an $S$-principal bundle over $X$.

\begin{proposition}\label{ample2quasiaff}
Let $\Lambda \subset \CDiv(X)$ be a polyhedral 
semigroup with associated toric bundle
$q \colon \t{X} \to X$ and ample locus $U \subset X$.
Then $q^{-1}(U) \cap \rq{X}$ is quasiaffine.
\end{proposition}

\begin{proof}
Consider the subgroup $\Gamma(\Lambda) \subset \CDiv(X)$ 
generated by $\Lambda$, 
and denote the associated graded 
$\mathcal{O}_{X}$-algebra by $\mathcal{B}$. 
Then $\rq{X}$ equals $\Spec(\mathcal{B})$,
and for any homogeneous $f \in \mathcal{B}(X)$,
its zero set as a function on $\rq{X}$ is 
equal to the inverse image $q^{-1}(Z(f)) \cap \rq{X}$.
Consequently, the set 
$q^{-1}(U) \cap \rq{X}$ is covered 
by affine open subsets of the form 
$\rq{X}_{f}$ with $f \in \mathcal{O}(\rq{X})$.
This gives the assertion. 
\end{proof}

We turn to the equivariant setting.
Let $G$ be a linear algebraic group, 
and suppose that $G$
acts on the normal prevariety $X$. 
A first observation is that the zero set $Z(f)$
of a homogeneous section $f$ behaves natural
with respect to the $G$-sheaf structure 
of Prop.~\ref{gsheaf}
arising from a $G$-linearization:

\begin{lemma}\label{translatezeroes}
Let $\Lambda \subset \WDiv(X)$ be a $G$-linearized 
polyhedral semigroup, and let $f \in \mathcal{A}_{D}(U)$
be a local section of the associated graded 
$\mathcal{O}_{X}$-algebra $\mathcal{A}$. Then we have
$Z(g \mal f) = g \mal Z(f)$ for any $g \in G$.
\end{lemma}

\begin{proof} 
By normality of $X$, we may assume that $U$ is smooth.
The problem being local, we may moreover assume
that $D$ is principal on $U$, say $D = -\Div(h)$.
Then the section $f$ is of the form $f = f'h$ with
a regular function $f'$,
and $Z(f)$ is just the zero set $Z(f')$ of $f'$.
Translating with $g \in G$ gives
$$ 
Z(g \mal f)
=
Z(g \mal f' \, g \mal h) 
=
Z(g \mal f') \cup Z(g \mal h).
$$

Since $h$ is a generator of $\mathcal{A}(U)$, 
the translate $g \mal h$ is a generator of 
$\mathcal{A}(g \mal U)$. 
This means that $Z(g \mal h)$ is empty. 
By definition of the $G$-sheaf structure,
$G$ acts canonically on the structure sheaf
$\mathcal{O}_{X}$, that means that $g \mal f' (x)$
equals $f'(g^{-1} \mal x)$.
This implies $Z(g \mal f') = g \mal Z(f')$,
and the assertion follows.
\end{proof}

\begin{proposition}\label{ampleinvar}
Let $\Lambda \subset \WDiv(X)$
be a $G$-linearized polyhedral semigroup. 
Then the Cartier locus and the ample locus of $\Lambda$
are $G$-invariant.
\end{proposition}

\begin{proof}
Let $\mathcal{A}$ denote the graded $\mathcal{O}_{X}$-algebra
corresponding to $\Lambda$.
The Cartier locus of $\Lambda$ is the set of all points
$x \in X$ such that for any $D \in \Lambda$ the stalk
$\mathcal{A}_{D,x}$ is generated by a single element.
Thus, using the $G$-sheaf structure of $\mathcal{A}$, 
we obtain that the Cartier locus is $G$-invariant.
Invariance of the ample locus, then is a simple 
consequence of Lemma~\ref{translatezeroes}.
\end{proof}

As a direct application, we extend a fundamental 
observation of Sumihiro on actions of connected 
linear algebraic groups $G$ on normal varieties 
$X$, see~\cite[Lemma~8]{Su}, and~\cite[Thm.~3.8]{Su2}: 
Every point $x \in X$ admits a $G$-invariant 
quasiprojective open neighbourhood.
Our methods give more generally:

\begin{proposition}\label{qptranslates}
Let $G$ be a connected linear algebraic group, let 
$X$ be a normal $G$-variety, and let $U \subset X$ be an 
open subset.
\begin{enumerate}
\item If $U$ is quasiprojective, then $G \mal U$ is 
   quasiprojective.
\item If $U$ is divisorial, then $G \mal U$ is 
   divisorial.
\end{enumerate}
In particular, the maximal quasiprojective and the maximal divisorial
open subsets of $X$ are $G$-invariant.
\end{proposition}

If $X$ admits a normal completion for which the factor
group of Weil divisors modulo $\QQ$-Cartier divisors is of finite
rank, then~\cite[Thm.~A]{Wl} says that $X$ has only finitely
many maximal open quasiprojective subvarieties. 
In particular, then~\ref{qptranslates}~(i) 
even holds with any connected algebraic group $G$, 
see~\cite[Thm.~D]{Wl}. 
A special case of the second statement is proved
in~\cite[Lemma~1.7]{AcHa}.

\begin{lemma}
\label{divext}
Let $X$ be a normal variety,
$D'$ a Weil divisor on some
open $U \subset X$,
and  $f'_{1}, \ldots, f'_{r}$ 
sections of $D'$ with 
$U \setminus Z(f'_{i})$ affine. 
Then there is a Weil divisor $D$ on $X$
allowing global sections $f_{1}, \ldots, f_{r}$ 
such that
$$
D\vert_{U} = D',
\qquad
f_{i} \vert_{U} = f'_{i},
\qquad
X \setminus Z(f_{i}) = U \setminus Z(f'_{i}).
$$
Moreover, if $U$ and $D'$ are invariant with respect 
to a given algebraic group action on $X$, then 
one can also choose $D$ to be so.
\end{lemma}

\begin{proof}
Let $D_{1}, \ldots, D_{s}$ be the prime divisors 
contained in $X \setminus U$.
Since the complement of $U  \setminus Z(f'_{i})$ 
in $X$ is of pure 
codimension one, we have
$$ 
U \setminus Z(f'_{i}) 
= 
X \setminus (D_{1} \cup \ldots \cup D_{s} \cup \b{Z(f'_{i})}).
$$
Consequently, by closing the components of $D'$ and adding
a suitably big multiple of $D_{1}+ \ldots + D_{s}$, 
we obtain the desired Weil divisor $D$ on $X$.
\end{proof}

\begin{proof}[Proof of Prop.~\ref{qptranslates}]
For (i), we choose a $D' \in \CDiv(U)$
allowing sections $f_{1}', \ldots, f_{r}'$
such that the sets $U \setminus Z(f_{i}')$ 
are affine and cover $U$.
Similarly, for (ii), we find
$D'_1, \ldots D'_r \in \CDiv(U)$
allowing sections $f_{1}', \ldots, f_{r}'$
such that the sets $U \setminus Z(f_{i}')$ 
are affine and cover $U$.

Use Lemma~\ref{divext} to extend $D'$ 
(resp. the $D_i'$)
to Weil divisors $D$ (resp. $D_i$) on $X$ such
that the $f_{i}'$ extend to global sections
$f_{i}$ over $X$ and satisfy 
$X \setminus Z(f_{i}) = U \setminus Z(f_{i}')$.
Let $\Lambda$ be the semigroup generated by $D$ 
(resp. the subgroup generated by the $D_i$).
Then, in both cases $U$ is contained in 
the ample locus of the extension $\Lambda$. 

Now, passing to subsemigroups of finite index does 
not shrink the ample locus.
Hence, we can use Prop.~\ref{existence1}, 
and endow $\Lambda$ with a $G$-linearization.
The assertion then follows from $G$-invariance of 
the ample locus of $\Lambda$ and the fact that 
quasiprojectivity as well as divisoriality transfer
to open subvarieties.
\end{proof}

We conclude this section with an
equivariant and refined version
of Prop.~\ref{ample2quasiaff}; again, we consider the 
subset $\rq{X} \subset \t{X}$ of free orbits
of the torus $S = \Spec(\KK[\Gamma(\Lambda)])$:

\begin{proposition}\label{ampleembed}
Let $\Lambda \subset \CDiv(X)$ be a $G$-linearized polyhedral
semigroup with associated toric bundle $q \colon \t{X} \to X$.
Let  $U \subset X$ be the ample locus of $\Lambda$,
and set $\rq{U} := \rq{X} \cap q^{-1}(U)$.
Then there is a $(G \times S)$-equivariant open embedding 
$\rq{U} \to Z$ into an affine $(G \times S)$-variety $Z$. 
Moreover, 
\begin{enumerate}
\item one can achieve that the image of the pullback map 
   $\mathcal{O}(Z) \to \mathcal{O}(\rq{U})$ 
   is contained in $\mathcal{O}(\t{X})$,
\item given $f_{1}, \ldots, f_{k} \in \mathcal{A}(X)$ as
  in~\ref{ampledef} with $\t{X}_{f_{i}} \subset \rq{X}$,
  one can achieve that each $f_{i}$ extends regularly to $Z$
  and satisfies $\rq{U}_{f_{i}} = Z_{f_{i}}$. 
\item For every $f \in \mathcal{O}(Z) \subset \mathcal{O}(\t{X})$ 
  with $\t{X}_f \subset \rq{X}$ and $f \vert_{Z \setminus \rq{U}} = 0$, 
  we have $Z_f = \rq{X}_f$.
\end{enumerate}
\end{proposition}

\begin{proof}
Let $f_{1}, \ldots, f_{k} \in \mathcal{A}(X)$ be 
as in~(ii), and complement this collection
by further homogeneous sections 
$f_{k+1}, \ldots, f_{r} \in \mathcal{A}(X)$
as in Def.~\ref{ampledef}
such that the affine sets $X_{i} := X \setminus Z(f_{i})$ 
cover the ample locus $U \subset X$.
Then each $f_{i}$, regarded as a regular function on $\t{X}$, 
vanishes outside the affine open set
$\t{X}_{i} := q^{-1}(X_{i})$ and has no zeroes inside
$\t{X}_{i} \cap \rq{X}$.

For each $i$, we choose finitely many homogeneous 
functions $h_{ij} \in \mathcal{O}(\t{X})$ such that
the affine algebra $\mathcal{O}(\t{X})_{f_{i}}$ is
generated by functions $h_{ij}/f_i^{l_{ij}}$.
Since the $G$-representation on $\mathcal{O}(\t{X})$
is rational, we find finite dimensional graded $G$-modules
$M_{i}, M_{ij} \subset \mathcal{O}(\t{X})$ such
that $f_{i} \in M_{i}$ and $h_{ij} \in M_{ij}$ holds.

Let $R \subset \mathcal{O}(\t{X})$ denote the subalgebra
generated by the elements of the $M_{i}$ and the $M_{ij}$.
Then $R$ is graded, $G$-invariant, and hence defines 
an affine $(G \times S)$-variety $Z := \Spec(R)$.
Note that $Z_{f_i} = \t{X}_{f_{i}}$ holds.
This gives $\rq{U} \subset Z$ and (ii).
Moreover, we obtain~(iii) by covering $\t{X}_f$
and $Z_f$ with the affine sets
$\t{X}_{f_if} = Z_{f_if}$.
\end{proof}

\section{Construction of quotients}

In this section, $G$ is a reductive group, 
and $X$ is a normal $G$-variety.
We describe the $G$-invariant 
open subsets $U \subset X$ admitting a good 
quotient with a quasiprojective or a divisorial
good quotient space.
First recall the precise definition of 
a good quotient, 
compare~\cite[p.~38]{Mu} and~\cite[Def.~1.5]{Se}:

\begin{definition}
A {\em good quotient\/} for a $G$-prevariety $X$ is an affine 
$G$-invariant morphism $p \colon X \to Y$ 
such that the canonical map 
$\mathcal{O}_{Y} \to p_{*}(\mathcal{O}_{X})^{G}$
is an isomorphism.
A good quotient is called {\em geometric\/} 
if its fibres are precisely the $G$-orbits.
\end{definition}

In our setting, a separated $G$-variety may have a 
good quotient with a nonseparated quotient space.
If a good quotient $X \to Y$ exists for a $G$-variety $X$,
then it is categorical, i.e., any $G$-invariant morphism 
$X \to Z$ factors uniquely through $X \to Y$. 
In particular, good quotient spaces are unique up to 
isomorphism. 
As usual, we write $X \to X \quot G$ for a good and 
$X \to X/G$ for a geometric quotient.

In general, the $G$-variety $X$ itself need not
admit a good quotient,
but there frequently exist many $G$-invariant open subsets 
$U \subset X$ with a good quotient.
Following~\cite{BB}, we say that 
a subset $V$ of an open $G$-invariant 
subset $U \subset X$ with good quotient
$p \colon U \to U \quot G$ is
{\em $G$-saturated\/} in $U$ if $V = p^{-1}(p(V))$
holds.

We begin with the construction of quasiprojective
good quotient spaces.
Fix a Weil divisor $D$ on $X$, and a $G$-linearization
of the semigroup $\Lambda := \NN D$;
we shall speak in the sequel of the $G$-linearized
Weil divisor $D$.
Recall from Prop.~\ref{gsheaf}
that there is an induced $G$-representation
on the global sections 
$\mathcal{A}(X)$ of the associated 
$\Lambda$-graded $\mathcal{O}_{X}$-algebra
$\mathcal{A}$.

\begin{definition}\label{semistabdef}
We call a point $x \in X$ {\em semistable\/} if there is 
an integer $n>0$ and a $G$-invariant 
$f \in \mathcal{A}_{nD}(X)$ such that 
$X \setminus Z(f)$ is a affine neighbourhood
of $x$ and $D$ is Cartier on $X \setminus Z(f)$.
\end{definition}

Following Mumford's notation, we denote the set of
semistable points of a $G$-linearized Weil divisor 
$D$ on $X$ by $X^{ss}(D)$, or by $X^{ss}(D,G)$ if 
we want to specify the group $G$.
Our concept of semistability yields all open subsets
admitting a quasiprojective good quotient space:  

\begin{theorem}\label{qpquots}
Let a reductive group $G$ act on a normal variety $X$.
\begin{enumerate}
 \item For any $G$-linearized Weil divisor $D$ on $X$, 
  there is a good quotient $X^{ss}(D) \to X^{ss}(D) \quot G$ 
  with a quasiprojective variety $X^{ss}(D) \quot G$.
\item If $U \subset X$ is open, $G$-invariant, and has 
   a good quotient $U \to U \quot G$ with $U \quot G$
   quasiprojective, then $U$ is a $G$-saturated subset of
   the set $X^{ss}(D)$ of semistable points of a canonically
   $G$-linearized Weil divisor $D$.
\end{enumerate}
\end{theorem}

\begin{proof} 
For i), we can follow the lines of~\cite[Thm.~1.10]{Mu}: 
Choose $G$-invariant homogeneous sections
$f_{1}, \ldots, f_{r} \in \mathcal{A}(X)$ as in
Def.~\ref{semistabdef} such that $X^{ss}(D)$ is covered 
by the sets $X_{i} := X \setminus Z(f_{i})$. 
Replacing the $f_{i}$ with suitable powers, we may assume that 
all of them have the same degree. 
Consider the good quotients:
$$ 
p_{i} \colon X_{i} \to X_{i} \quot G 
= \Spec(\mathcal{O}(X_{i})^{G}).
$$

Each $X_{i} \setminus X_{j}$ is the zero set of the $G$-invariant 
regular function $f_{j}/f_{i}$. 
Thus $X_{i} \cap X_{j}$ is saturated
with respect to the quotient map 
$p_{i} \colon X_{i} \to X_{i} \quot G$. 
It follows that the $p_{i}$ glue 
together to a good quotient  
$p \colon X^{ss}(D) \to X^{ss}(D) \quot G$.
Moreover, for fixed $i_0$, the $f_{i_{0}}/f_{i}$ are local 
equations for an ample divisor on $X^{ss}(D) \quot G$.

To prove~(ii), let $Y := U \quot G$, and let $p \colon U \to Y$ 
be the quotient map. 
Choose an ample divisor $E$ on $Y$ allowing
global sections $h_{1}, \ldots, h_{r}$ such that 
the sets $Y \setminus Z(h_{i})$ form an affine cover 
of $Y$.
Consider the pullback data
$D' := p^{*}E$ and $f'_{i} := p^{*}(h_{i})$.
Then Lemma~\ref{divext} provides a $G$-invariant 
Weil divisor $D$ on $X$ extending $D'$ and
sections $f_{i}$ extending $f'_{i}$
such that
$$ 
X \setminus Z(f_{i}) 
\; = \;
U \setminus Z(f'_{i})
\; = \; 
p^{-1}(Y \setminus Z(h_{i})).
$$

Let $\mathcal{A}$ be the graded $\mathcal{O}_{X}$-algebra
associated to $D$, and consider the canonical $G$-linearization of $D$ 
provided by Prop.~\ref{existence2}.
Then the sections $f_{i} \in \mathcal{A}(X)$ are $G$-invariant,
and satisfy the conditions of~\ref{semistabdef}.
It follows that $U$ is a saturated subset of $X^{ss}(D)$.
\end{proof}

For the construction of divisorial quotient spaces,
we work with finitely generated subgroups
$\Lambda \subset \WDiv(X)$; 
these are in particular polyhedral semigroups.
Fix such a subgroup $\Lambda \subset \WDiv(X)$, 
and a $G$-linearization of $\Lambda$ as introduced
in Section~\ref{section2}.
Again, we have an induced $G$-representation 
on the global sections $\mathcal{A}(X)$ of 
the associated $\Lambda$-graded $\mathcal{O}_{X}$-algebra  
$\mathcal{A}$.

\begin{definition}\label{semistabgroupdef}
We call a point $x \in X$ {\em semistable\/}, 
if $x$ has an affine neighbourhood 
$U = X \setminus Z(f)$ with some $G$-invariant 
homogeneous $f \in \mathcal{A}(X)$
such that all $D \in \Lambda$ are Cartier on $U$,
and the $D \in \Lambda$ admitting a $G$-invariant
invertible $h \in \mathcal{A}_{D}(U)$
form a subgroup of finite index in $\Lambda$.
\end{definition}

As before, the set of semistable points is denoted by 
$X^{ss}(\Lambda)$, or $X^{ss}(\Lambda,G)$
if we want to specify the group $G$. 
Note that for $G$-linearized groups of Cartier divisors
we retrieve the notion of semistability 
introduced in~\cite[Def.~2.1]{Ha0}.
We obtain the following generalizations 
of~\cite[Thms.~3.1, 4.1]{Ha0}:

\begin{theorem}\label{divisorialquotients}
Let a reductive group $G$ act on a normal variety $X$.
\begin{enumerate}
 \item For any $G$-linearized group $\Lambda \subset \WDiv(X)$,
  there is a good quotient $X^{ss}(\Lambda) \to X^{ss}(\Lambda) \quot G$ 
  with a divisorial prevariety $X^{ss}(D) \quot G$.
\item If $U \subset X$ is open, $G$-invariant, and admits 
  a good quotient $U \to U \quot G$ 
  with $U \quot G$ divisorial, then $U$ is a $G$-saturated
  subset of the set $X^{ss}(\Lambda)$ of semistable points of a
  canonically $G$-linearized group $\Lambda \subset \WDiv(X)$.
\end{enumerate}
\end{theorem}

\begin{proof} To prove (i), consider the Cartier locus 
$X_{0} \subset X$ of $\Lambda$.
By Prop.~\ref{ampleinvar}, the set $X_{0}$ is $G$-invariant.
Since $X$ is normal, $X \setminus X_{0}$ is of codimension
at least two in $X$.
Hence $X_{0}^{ss}(\Lambda)$ equals $X^{ss}(\Lambda)$,
and we may assume that $\Lambda$ consists of Cartier
divisors.
But then~\cite[Thm.~3.1]{Ha0} gives the assertion.

The proof of (ii) is analogous to that of~\cite[Thm.~4.1]{Ha0}.
Using divisoriality of $Y := U \quot G$ 
and~\cite[Lemma~4.3]{Ha0}, 
we find effective $E_{1}, \ldots, E_{r} \in \CDiv(Y)$ 
and global sections $h_{ij}$ of the $E_{i}$ such that
the sets $V_{ij} := Y \setminus Z(h_{ij})$ 
form an affine cover of $Y$, 
and every $E_{k}$ admits an invertible section $h_{ijk}$ 
over $V_{ij}$. 

Let $p \colon U \to Y$ be the quotient map.
Lemma~\ref{divext} provides invariant Weil
divisors $D_{i}$ on $X$ admitting global sections
$f_{ij}$ such that with $U_{ij} := p^{-1}(V_{ij})$
we have 
$$ 
D_{i} \vert_{U}
\; = \;
p^{*}E_{i},
\qquad
f_{ij}\vert_{U_{ij}} 
\; = \; 
p^{*}(h_{ij}),
\qquad 
X \setminus Z(f_{ij}) 
\; = \;
p^{-1}(V_{ij}). 
$$ 
By Prop.~\ref{existence2}, 
the group $\Lambda \subset \WDiv(X)$ 
generated by the $D_{i}$ is canonically
$G$-linearized.
The sections $f_{ij}$ and $p^{*}(h_{ijk})$
serve to verify $U \subset X^{ss}(\Lambda)$.
Since each $U_{ij}$
is $G$-saturated in $X^{ss}(\Lambda)$,
the same holds for $U$.
\end{proof} 

We conclude the section with an example,
showing that in the singular case
Mumford's method and the generalization
given in~\cite{Ha0} need no longer provide all open subsets
with quasiprojective or divisorial quotient spaces.
Consider the cone $X$ over the image of 
$\PP_{1} \times \PP_{1}$ in $\PP_{3}$ 
under the Segre embedding, i.e.:
$$ X = V(\KK^{4}; z_{1}z_{3} - z_{2}z_{4}). $$

Then $X$ is a normal variety having precisely one singular point.
Let $U := X_{z_{2}} \cup X_{z_{4}}$ be the set of points having
nonvanishing 2nd or 4th coordinate.
We consider the following action of the twodimensional torus 
$T := \KK^{*} \times \KK^{*}$ on $X$:
$$
t \mal x := 
(t_{1}^{2}x_{1},
t_{1}t_{2}^{2}x_{2},
t_{1}t_{2}x_{3},
t_{1}^{2}t_{2}^{-1}x_{4}).
$$

\begin{proposition}\label{counterexample}
The set $U \subset X$ has a geometric quotient $U \to U / T$ 
with $U / T \cong \PP_{1}$, but $U$ is not the 
set of semistable points of a $T$-linearized line 
bundle on $X$ in the sense of~\cite[Def.~1.7]{Mu}.
\end{proposition}

\begin{proof}
The most convenient way is to view $X$ as a toric variety, 
and to work in the language of lattice fans, see~\cite{Fu} 
for the basic notions. 
As a toric variety, $X$ corresponds to the lattice cone $\sigma$ 
in $\ZZ^{3}$ generated by the vectors
$$ 
v_{1} := (1,0,0),
\qquad
v_{2} := (0,1,0),
\qquad
v_{3} := (0,1,1),
\qquad
v_{4} := (1,0,1).
$$

The big torus of $X$ is $T_{X} = (\KK^{*})^{3}$.
The torus $T$ acts on $X$ by $(t,x) \mapsto \varphi(t) \mal x$,
where $\varphi \colon T \to T_{X}$ is the homomorphism of tori
corresponding to the linear map
$$ 
\ZZ^{2} \to \ZZ^{3},
\qquad
(1,0) \mapsto (2,1,1),
\qquad
(0,1) \mapsto (0,2,1).
$$

Our open set $U \subset X$ is a union of three
$T_{X}$-orbits: the big $T_{X}$-orbit, and the 
two twodimensional $T_{X}$-orbits corresponding to the rays 
$\varrho_{1} := \QQ_{\ge 0} v_{1}$ and 
$\varrho_{3} := \QQ_{\ge 0} v_{3}$ of the cone $\sigma$.
The fan theoretical criterion~\cite[Thm.~5.1]{Hm},
tells us that there is a geometric quotient for the action 
of $T$ on $U$; 
namely the toric morphism $p \colon U \to \PP_{1}$ 
defined by the linear map
$$ 
P \colon \ZZ^{3} \to \ZZ,
\qquad
(w_{1},w_{2},w_{3})
\mapsto 
w_{1} + 2w_{2} - 4w_{3}.
$$

We show now that there is no $T$-linearized line bundle on $X$ having 
$U$ as its set of semistable points. First note that as an affine
toric variety, $X$ has trivial Picard group. Thus we only have to
consider $T$-linearizations of the trivial bundle. Since
$\mathcal{O}^{*}(X) = \KK^{*}$ holds, each such linearization
is given by a character $\chi$ of $T$:
$$ t \mal (x,z) = (t \mal x, \chi(t)z). $$

Consequently, in view of~\cite[Def.~1.7]{Mu},
we have to show that $U$ is not a union of sets $X_{f}$,
for a collection of functions $f \in \mathcal{O}(X)$ that are
$T$-homogeneous with respect to a common character of 
the torus $T$.

Now, any $T$-homogeneous regular function on $X$
is a sum of $T$-homo\-ge\-neous character functions 
$\chi^{u} \in \mathcal{O}(X)$ , 
where $u = (u_{1},u_{2},u_{3})$ is a lattice vector of the dual 
cone $\sigma^{\vee}$ of $\sigma$. 
Recall that $u \in \sigma^{\vee}$ means that the linear form $u$
is nonnegative on $\sigma$, i.e., we have
$$ 
u_{1} \ge 0,
\qquad
u_{2} \ge 0,
\qquad
u_{2} + u_{3} \ge 0,
\qquad
u_{1} + u_{3} \ge 0.
$$

For such a character function $\chi^{u} \in \mathcal{O}(X)$, we can
determine its weight with respect to $T$ by applying the dual of 
the embedding $\ZZ^{2} \to \ZZ^{3}$ to the vector $u$. 
Thus, $\chi^{u}$ is $T$-homogeneous with respect to the character 
of $T$ corresponding to the lattice vector
$$ (2u_{1}+u_{2}+u_{3},2u_{2}+u_{3}). $$

The conditions that a character function 
$\chi^{u} \in \mathcal{O}(X)$ does not vanish along 
the orbit $T_{X} \mal x_{i}$ corresponding 
to one of the rays $\varrho_{i}$ are $u_{1}= 0$ for nonvanishing along
$T_{X} \mal x_{1}$, and $u_{3} = -u_{2}$ for nonvanishing along $T_{X}
\mal x_{3}$. 

Suppose that $\chi^{u}\in \mathcal{O}(X) $ does not vanish along 
$T_{X} \mal x_{1}$ and that $\chi^{\t{u}} \in \mathcal{O}(X)$ 
does not vanish along $T_{X} \mal x_{3}$.
Then their respective $T$-weights are given by the vectors
$$
(u_{2}+u_{3},2u_{2}+u_{3}),
\qquad
(2 \t{u}_{1}, \t{u}_{2}).
$$
If both are $T$-homogeneous with respect to the same character, 
then we must have $2\t{u}_{1} \le \t{u}_{2}$. 
But then nonvanishing along $T_{X} \mal x_{3}$ and the last regularity
condition imply $\t{u} = 0$.

In conclusion, we obtain that only the trivial character of $T$
admits homogeneous functions that do not vanish along 
$T_{X} \mal x_{1}$ and functions that do not vanish along 
$T_{X} \mal x_{3}$. 
Since $T$ acts with an attractive fixed point on $X$ this means
that we cannot obtain $U$ as a union of sets $X_{f}$ as needed.
\end{proof}

\section{First Hilbert-Mumford type statements}

We come to the first Hilbert-Mumford type result of the article.
It allows us to express the set of $G$-semistable points in terms of 
the $T$-semistable points for a maximal torus $T \subset G$.
In the case of an ample divisor $D$ on a projective $G$-variety, 
the first assertion of our result is equivalent to~\cite[Thm.~2.1]{Mu}.

\begin{theorem}\label{himusemistab}
Let a reductive group $G$ act on a normal 
variety $X$, and let $T \subset G$ be a maximal torus.
\begin{enumerate} 
\item Let $D$ be a $G$-linearized Weil 
divisor on $X$. Then we have:
$$ X^{ss}(D,G) = \bigcap_{g \in G} g \mal X^{ss}(D,T). $$
\item Let 
$\Lambda \subset \WDiv(X)$ be a $G$-linearized 
subgroup. Then we have:
$$ 
X^{ss}(\Lambda,G) 
= 
\bigcap_{g \in G} g \mal X^{ss}(\Lambda,T). 
$$
\end{enumerate}
\end{theorem}

The proof relies on a geometric analysis of instability;
it makes repeated use of the classical Hilbert-Mumford
Theorem, see for example~\cite[Thm.~4.2]{Bi}:

\begin{theorem}
\label{HiMuclass}
Let a reductive group $G$ act on an affine 
variety $Z$, let $z \in Z$, and let $Y \subset \b{G \mal z}$ 
be a $G$-invariant closed subset.
Then there is a one parameter 
subgroup $\lambda \colon \KK^* \to G$ with
$\lim_{t \to 0} \lambda(t) \mal z \in Y$.
\end{theorem}

The basic preparing steps concern
the following situation:
$G$ is a reductive group, 
$Z$ is an affine $G$-variety,
and $T \subset G$ is a maximal torus.
Then we have good quotients
$$ 
p_{T} \colon Z \to Z \quot T,
\qquad 
p_{G} \colon Z \to Z \quot G. 
$$

\begin{lemma}\label{instab1}
Let $A \subset Z$ be $G$-invariant and closed,
and let $z \in p_{G}^{-1}(p_{G}(A))$.
Then there is a $g \in G$ with
$g \mal z \in p_{T}^{-1}(p_{T}(A))$.
\end{lemma}

\begin{proof}
Since $p_{G} \colon Z \to Z \quot G$
separates disjoint $G$-invariant closed sets,
the closure of $G \mal z$ intersects $A$. 
By Theorem~\ref{HiMuclass},
there is a maximal torus $S \subset G$ such that
the closure of $S \mal z$ intersects $A$.
Choose a $g \in G$ with $T = gSg^{-1}$.
Then the closure of $T \mal g \mal z$
intersects $A$.
This implies 
$p_{T}(g \mal z) \in p_{T}^{-1}(p_{T}(A))$.
\end{proof}

Suppose that in addition to the $G$-action
there is an action of $\KK^{*}$ on $Z$
such that these two actions commute.
Then there are induced $\KK^{*}$-actions on the quotient
spaces $Z \quot T$ and $Z \quot G$
making the respective quotient maps equivariant.
Let $B_{T}^{0} \subset Z \quot T$ and 
$B_{G}^{0} \subset Z \quot G$ denote 
the fixed point sets of these $\KK^{*}$-actions.

\begin{lemma}\label{instab2}
Let $z \in Z$ with $p_{G}(z) \in B_{G}^{0}$.
Then there is a $g \in G$ with
$p_{T}(g \mal z) \in B_{T}^{0}$.
\end{lemma}

\begin{proof} 
Let $G \mal z_{0}$ be the closed $G$-orbit
in the fibre $p_{G}^{-1}(p_{G}(z))$. 
If $z_{0}$ is a fixed point of the $\KK^{*}$-action
on $Z$, then the whole orbit $G \mal z_{0}$ consists 
of $\KK^{*}$-fixed points, and the assertion is a direct 
consequence of Theorem~\ref{HiMuclass}.
So we may assume for this proof 
that the orbit $\KK^{*} \mal z_{0}$ is 
nontrivial. 

By Theorem~\ref{HiMuclass}, there is a onedimensional 
subtorus $S_0 \subset G$ and a $g_0 \in G$ such that $z_0$ 
lies in the closure of $S_0 \mal z'$, where 
$z' := g_0 \mal z$.
Note that for any $t \in \KK^*$, 
the point $t \mal z_{0}$ lies in the closure of 
$S_{0} \mal t \mal z'$.
This implies in particular that any point of 
$\KK^{*} \mal z_{0}$ is fixed by $S_{0}$.
Consequently, $S_{0}$ is a subgroup of the 
stabilizer $G_{0}$ of $\KK^{*} \mal z_{0}$.

Let $n \in \NN$ denote the order of the
isotropy group of $\KK^{*}$ in $z_{0}$.
Then the orbit maps
$\mu \colon g \mapsto g \mal z_0$ of $G_0$ 
and $\nu \colon t \mapsto t \mal z_0$
of $\KK^{*}$ give rise to a well defined
morphism of linear algebraic groups:
$$ 
G_{0} \to \KK^{*},
\qquad
g \mapsto 
(\nu^{-1}(\mu(g)))^{n}.
$$

Clearly, $S_{0}$ is contained in the kernel of this
homomorphism.
By general properties of linear algebraic groups
any maximal torus of $G_{0}$ is mapped onto $\KK^{*}$,
see e.g.~\cite[Cor.~C, p.~136]{Hu}.
We choose a maximal torus $S_{1} \subset G_0$ such
that $S_1$ contains~$S_{0}$.

Let $S \subset G$ be a maximal torus with
$S_1 \subset S$.
Then $z_{0}$ lies in the closure of $S \mal z'$.
Moreover, $\KK^{*} \mal z_{0}$ is contained in
$S \mal z_{0}$.
Writing $S = g_1^{-1}Tg_1$ with a suitable $g_{1} \in G$, 
we obtain that $g_1 \mal z_0$ lies in the closure of
$T \mal g_1 \mal z' $, and $\KK^{*} \mal g_1 \mal z_0$
is contained in $T \mal g_1 \mal z_{0}$.
Thus, $g := g_1g_0$ is as wanted.
\end{proof}

The next observation concerns limits with respect to the 
$\KK^{*}$-action on the quotient spaces.
For $H = T$ and $H = G$ we consider the sets:
$$ 
B_{H}^{-} 
:= 
\{y \in Z \quot H; \; \lim_{t \to \infty} t \mal y
           \text{ exists and differs from } y\}. 
$$

\begin{lemma}\label{instab3}
Let $z \in Z$ with $p_{G}(z) \in B_{G}^{-}$.
Then there is a $g \in G$ such that 
$p_{T}(g \mal z) \in B_{T}^{-}$ holds.
\end{lemma}

\begin{proof}
Let $y_{0} \in Z \quot G$ be the
limit point of $p_{G}(z)$,
and choose $z_{0} \in Z$ with
$G \mal z_{0}$ closed in $Z$ and $p_{G}(z_{0}) = y_{0}$.
Note that $G \mal z_{0}$ is $\KK^{*}$-invariant.
Consider the quotient $q \colon Z \to Z \quot \KK^{*}$.
Then $G \mal q(z_{0})$ is contained
in the closure of $G \mal q(z)$,
because $q(G \mal z_{0})$ is closed, and we have
$$  (Z \quot G) \quot \KK^{*} = (Z\quot \KK^{*}) \quot G. $$

Thus, according to Theorem~\ref{HiMuclass}, there
exist $g, g_0 \in G$ such that 
$g_0 \mal q(z_{0})$ lies in the closure of
$T \mal g \mal q(z)$.
We can conclude that in $Z \quot T$,
the $\KK^{*}$-orbit closures of the points $p_{T}(g \mal z)$ 
and $p_{T}(g_{0} \mal z_{0})$ intersect nontrivially; 
this time we use
$$  (Z \quot T) \quot \KK^{*} = (Z\quot \KK^{*}) \quot T. $$

Since we have a $\KK^{*}$-equivariant map
$Z \quot T \to Z \quot G$,
and there is a $G$-invariant homogeneous 
function $f \in \mathcal{O}(Z)$ of negative
weight with $f(z) \ne 0$ and $f(z_{0}) = 0$,
it follows that $p_{T}(g \mal z)$ belongs to
$B_{T}^{-}$.
\end{proof} 

\begin{proof}[Proof of Theorem~\ref{himusemistab}]
For (i), we may assume that $D$ is nontrivial.
By Prop.~\ref{ampleinvar},
the Cartier locus $X_{0} \subset X$ 
of $\Lambda := \NN D$ is $G$-invariant.
Moreover, by normality of $X$, the complement
$X \setminus X_{0}$ is of codimension at least 
two in $X$.
Consequently, $X_{0}^{ss}(D,T)$ equals $X^{ss}(D,T)$,
and $X_{0}^{ss}(D,G)$ equals $X^{ss}(D,G)$.
Thus we may assume for this proof that $X = X_{0}$
holds.

Let $\mathcal{A}$ be the graded $\mathcal{O}_{X}$-algebra 
associated to $\Lambda$.
The associated $\t{X} := \Spec(\mathcal{A})$
is a line bundle over $X$, 
and the torus acting on $\t{X}$ is $\KK^{*}$.
Consider the $G$-action on 
$\t{X}$ provided by Prop.~\ref{liftedaction}.
Removing the zero section gives 
the $(G\times \KK^{*})$-invariant open 
subvariety $\rq{X} \subset \t{X}$.
Let $q \colon \t{X} \to X$ be the canonical 
$G$-equivariant map,
$U \subset X$ the ample locus of $D$,
and $\rq{U} := q^{-1}(U) \cap \rq{X}$.
 
Choose $T$-invariant homogeneous 
$f_{1}, \ldots, f_{r} \in \mathcal{A}(X)$ 
and $G$-invariant homogeneous 
$h_{1}, \ldots, h_{s} \in \mathcal{A}(X)$ 
as in Def.~\ref{semistabdef} 
such that the sets 
$X \setminus Z(f_{i})$ and $X \setminus Z(h_{j})$ 
cover $X^{ss}(D,T)$ and 
$X^{ss}(D,G)$ respectively.
Regarded as functions on $\t{X}$,
the $f_{i}$ and the $h_{j}$ vanish 
along the zero section $\t{X} \setminus \rq{X}$,
because they are of positive degree.

According to Prop.~\ref{ampleembed},
we can choose an equivariant 
open embedding  $\rq{U} \subset Z$ into an affine 
$(G \times \KK^*)$-variety
$Z$ with the following two properties:
Firstly, we have 
$\mathcal{O}(Z) \subset \mathcal{O}(\t{X})$.
Secondly, the functions $f_{i}, h_{j} \in \mathcal{O}(\rq{U})$
extend regularly to $Z$ and satisfy
$\rq{U}_{f_{i}} = Z_{f_{i}}$ and $\rq{U}_{h_{j}} = Z_{h_{j}}$. 

Now consider the induced $\KK^{*}$-actions on the quotient
spaces $Z \quot T$ and $Z \quot G$.
As before, let $B_{T}^{0}$, $B_{G}^{0}$ be the fixed point sets
of these $\KK^{*}$-actions,
and let $B_{T}^{-}$, $B_{G}^{-}$ be the sets of non fixed points
admitting a limit for $t \to \infty$.
Then, setting $A := Z \setminus \rq{U}$,
we claim that for the respective sets
of semistable points one has:
\begin{eqnarray*}
\rq{X} \cap q^{-1}(X^{ss}(D,T)) 
& = & 
Z \setminus p_{T}^{-1}(p_{T}(A)
               \cup B_{T}^{0}
               \cup B_{T}^{-}), \\
\rq{X} \cap q^{-1}(X^{ss}(D,G)) 
& = & 
Z \setminus p_{G}^{-1}(p_{G}(A)
               \cup B_{G}^{0}
               \cup B_{G}^{-}).
\end{eqnarray*}

Indeed, the inclusion ``$\subset$'' of the first equation
is due to the facts that the intersection 
$\rq{X} \cap q^{-1}(X \setminus Z(f_{i}))$
equals $Z_{f_{i}}$, and that each $f_{i}$ by $T$-invariance and
homogeneity of positive degree vanishes along 
the set $p_{T}^{-1}(p_{T}(A) \cup B_{T}^{0} \cup B_{T}^{-})$.
Analogously one obtains the inclusion ``$\subset$'' for the second 
equation. 

To see the inclusions ``$\supset$'', we treat again 
exemplarily the first equation.
The ideal of 
$p_{T}(A) \cup B_{T}^{0} \cup B_{T}^{-}$ in 
$\mathcal{O}(Z \quot T)$
is generated by functions $f'$ that are homogeneous 
of positive degree.
Since $\mathcal{O}(Z) \subset \mathcal{O}(\t{X})$ holds,
each $f := p_{T}^{*}(f')$ is a $T$-invariant homogeneous 
section of positive degree in $\mathcal{A}(X)$.
By Prop.~\ref{ampleembed}, we have
$$
Z_{f}
\; = \;
\rq{X}_{f}
\; = \;
\rq{X} \cap q^{-1}(X \setminus Z(f)).
$$

It follows that $X \setminus Z(f)$ is affine, and hence $f$
is as in Def.~\ref{semistabdef}.
Consequently, $Z_{f}$ lies over the set of $T$-semistable points 
of $X$. Since the functions $f$ generate the ideal of 
$p_{T}^{-1}(p_{T}(A) \cup B_{T}^{0} \cup B_{T}^{-})$,
we obtain the desired inclusion.

Now, Lemmas~\ref{instab1}, \ref{instab2}, and~\ref{instab3}
show that the inclusion ``$\supset$'' of the assertion is 
valid. The reverse inclusion is easy: Every translate
$g \mal X^{ss}(D,T)$ is the set of semistable points 
of $g T g^{-1}$ and hence contains $X^{ss}(D,G)$.

The proof of (ii) is similar.
As in the proof of~(i),
we may assume that $\Lambda$ consists 
of Cartier divisors. 
Let $\mathcal{A}$ be the associated $\Lambda$-graded 
$\mathcal{O}_{X}$-algebra.
Consider $\rq{X} := \Spec(\mathcal{A})$ 
with its actions of $S := \Spec(\KK[\Lambda])$ and $G$,
and the $G$-equivariant canonical map $q \colon \rq{X} \to X$.
Let $U \subset X$ be the ample locus of $\Lambda$,
and set $\rq{U} := q^{-1}(U)$.

Cover $X^{ss}(\Lambda,T)$ by 
sets $X \setminus Z(f_i)$ with
$T$-invariant homogeneous
$f_{i} \in \mathcal{A}(X)$
as in Def.~\ref{semistabgroupdef}. 
Similarly, cover $X^{ss}(\Lambda,G)$
by $X \setminus Z(h_j)$ 
with $G$-invariant homogeneous 
$h_{j} \in \mathcal{A}(X)$.
Lemma~\ref{ampleembed} provides
an equivariant open embedding $\rq{U} \subset Z$ 
into an affine $(G\times S)$-variety $Z$ 
with $\mathcal{O}(Z) \subset \mathcal{O}(\rq{X})$
such that all $f_{i}, h_{j}$
extend regularly to $Z$, 
and satisfy
$
\rq{U}_{f_{i}}
=
Z_{f_{i}}
$
and 
$
\rq{U}_{h_{j}}
=
Z_{h_{j}}.
$ 

For $H = T, G$, 
consider the quotient $p_{H} \colon Z \to Z \quot H$,
and induced action of $S$ on $Z \quot H$.
We describe $X^{ss}(\Lambda,H)$
in terms of these data. 
Let $A := Z \setminus \rq{U}$,
and let $B_{H}^{0} \subset Z \quot H$ 
be the set of all $y \in Z \quot H$
with an infinite isotropy group $S_{y}$.
We claim
\begin{equation}
\label{proofclaim2}
q^{-1}(X^{ss}(\Lambda,H)) 
\; = \; 
Z \setminus p_{H}^{-1}(p_{H}(A) \cup B_{H}^{0}). 
\end{equation}

The inclusion ``$\subset$'' is~\cite[Prop.~2.3~(i)]{Ha0}.
For the reverse inclusion, we use~\cite[Lemma~2.4]{Ha0}:
It tells us that the ideal of $p_{H}(A) \cup B_{H}^{0}$
in $\mathcal{O}(Z \quot H)$ is generated by $S$-homogeneous
elements $f'$ such that $\mathcal{O}(Z \quot H)_{f'}$ admits 
homogeneous invertible elements for almost every character
of the torus $S$.

For such $f'$, the pullback $f := p_{H}^{*}(f')$ 
is an $H$-invariant element of $\mathcal{O}(Z)$ and hence
of $\mathcal{A}(X)$, and, by Lemma~\ref{ampleembed}~(iii),
we have $\rq{X}_{f} = Z_{f}$.
Thus $q(\rq{X}_{f}) = X \setminus Z(f)$ is affine, 
and we see that $f$ is as in Def.~\ref{semistabgroupdef}.
Hence, $q^{-1}(X^{ss}(\Lambda,H))  \supset Z_{f}$
holds, which finally gives the claim.

Now, $B_{H}^{0} \subset  Z \quot H$ 
is the union of the fixed point sets 
$B_{H}^{0}(\mu)$ of all one parameter 
subgroups $\mu \colon \KK^{*} \to S$.
Lemmas~\ref{instab1} and~\ref{instab2}
tell us 
$$
p_G^{-1}(p_G(A)) 
\; = \; 
\bigcup_{g \in G} g \mal p_T^{-1}(p_T(A)),
\qquad
p_G^{-1}(B_G^{\circ}(\mu)) 
\; = \; 
\bigcup_{g \in G} g \mal p_T^{-1}(B_T^{\circ}(\mu)).
$$
Together with~(\ref{proofclaim2}), this gives
``$\supset$'' in the assertion.
The reverse inclusion is due to the fact 
that $g \mal X^{ss}(\Lambda,T)$ is the set of 
semistable points of $g T g^{-1} \subset G$.
\end{proof}

\section{Actions of semisimple groups}\label{section7}

In this section, we apply our results to
actions of semisimple groups.
This gives generalizations of several results 
presented in~\cite{BBSw2}, \cite{BBSw3}, and~\cite{Ha2}.
We work with
the following notions of maximality,
compare~\cite{BBSw2} and~\cite{Ha2}:

\begin{definition}
Let $G$ be a reductive group, let $X$ be a $G$-variety,
and let $U \subset X$ be a $G$-invariant open subset.
We say that
\begin{enumerate}
\item $U$ is a {\em qp-maximal $G$-set}
   if there is a good quotient $U \to U \quot G$ with
   $U \quot G$ quasiprojective, and $U$ is not a $G$-saturated
   subset of a properly larger $U' \subset X$ admitting a good
   quotient $U' \to U' \quot G$ with $U' \quot G$ quasiprojective,
\item $U$ is a {\em d-maximal $G$-set}
   if there is a good quotient $U \to U \quot G$ with
   $U \quot G$ divisorial, and $U$ is not a $G$-saturated
   subset of a properly larger $U' \subset X$ admitting a good
   quotient $U' \to U' \quot G$ with $U' \quot G$ divisorial.
\end{enumerate}
\end{definition}

In the sequel, $G$ is a connected semisimple group,
$T \subset G$ a maximal torus, and
$N \subset G$ the normalizer of $T$ in $G$.
Moreover, $X$ is a normal $G$-variety.
The first result is a further
Hilbert-Mumford type statement.
It generalizes~\cite[Cor.~1]{BBSw3}, and the 
results on the case $G = \SL_{2}$ 
given in~\cite[Thm.~9]{BBSw2} and~\cite[Thm.~2.2]{Ha2}:

\begin{theorem}\label{semisimhimu}
Let $U \subset X$ be an $N$-invariant open subset 
of $X$, and let $W(U)$ denote the intersection of all
translates $g \mal U$, where $g \in G$.
\begin{enumerate}
\item If $U \subset X$ is a qp-maximal $T$-set,
   then $W(U)$ is open and $T$-saturated in $U$,
   and there is a good quotient $W(U) \to W(U) \quot G$
   with $W(U) \quot G$ quasiprojective.
\item If $U \subset X$ is a d-maximal $N$-set,
   then $W(U)$ is open and $T$-saturated in $U$,
   and there is a good quotient $W(U) \to W(U) \quot G$
   with  $W(U) \quot G$ divisorial.
\end{enumerate}
\end{theorem}

The proof of this Theorem consists of combining the Hilbert-Mumford 
Theorems~\ref{himusemistab} with the following observation:

\begin{proposition}\label{max2semistab}
Let $U \subset X$ be an $N$-invariant open subset.
\begin{enumerate}
\item If $U$ is a qp-maximal $N$-set, then there 
   exists a $G$-linearized Weil divisor $D$ on 
   $X$ with $U = X^{ss}(D,N)$.
\item If $U$ is a d-maximal $N$-set, then there 
   is a $G$-linearized group $\Lambda \subset \WDiv(X)$
   with $U = X^{ss}(\Lambda,N)$.
\end{enumerate}
\end{proposition}

\begin{proof}
We exemplarily prove the first assertion.
By Theorem~\ref{qpquots}~(ii), 
there is a canonically $N$-linearized 
Weil divisor $D$ on $X$ such that $U$ 
is $N$-saturated in $X^{ss}(D,N)$. 
By qp-maximality of $U$, this implies
$U = X^{ss}(D,N)$. 
We show now that after possibly replacing $D$ with a positive 
multiple, the $N$-lineari\-za\-tion extends to a 
$G$-linearization.

Let $Z$ be a $G$-equivariant 
completion of $X$, see~\cite[Theorem~3]{Su}. 
Applying equivariant normalization, we achieve that
$Z$ is normal. 
By closing the support, we extend $D$ 
to a Weil divisor $E$ of $Z$.
Then $E$ is $N$-invariant and hence, by 
Prop.~\ref{existence2}, 
it is canonically $N$-linearized.

Prop.~\ref{existence1} tells us that after
replacing $E$ (and $D$) with a suitable multiple, 
we can choose a $G$-linearization of $E$.
Since we have $\mathcal{O}(Z) = \KK$ 
and the character group of $N$ is finite,
Prop.~\ref{uniqueness}~(ii)
says that after possibly passing to a further multiple,
the $G$-linearization of $E$ induces 
the canonical $N$-linearization of $E$ over $Z$.
Restricting to $X \subset Z$, we obtain the assertion.
\end{proof}

Note that this proposition is the place, 
where semisimplicity of $G$ came in.
In the proof, we made essentially use of the 
fact that the character group of $N$ is finite.

\begin{proof}[Proof of Theorem~\ref{semisimhimu}]
Note first that in the setting of (i), the induced action of 
the Weyl group $N/T$ on $U \quot T$ admits a 
geometric quotient with a quasiprojective quotient 
space.
The composition of the quotients by $T$ and $N/T$
is a good quotient $U \to U \quot N$.
It follows that $U$ is a qp-maximal $N$-set.

Now, for~(i), choose a $G$-linearized semigroup $\Lambda= \NN D$,
and, for ~(ii), a $G$-linearized group $\Lambda \subset \WDiv(X)$ 
as provided by Prop.~\ref{max2semistab}.
By the definition of semistability, we have
$$ 
X^{ss}(\Lambda,G)
\; \subset \;
\bigcap_{g \in G} g \mal X^{ss}(\Lambda,N)
\; \subset \;
\bigcap_{g \in G} g \mal X^{ss}(\Lambda,T).
$$
{}From Theorem~\ref{himusemistab} 
we infer that also
the reverse inclusions hold.
This gives the assertion.
\end{proof}

In the case of complete quotient spaces,
the approach via Weil divisors finally turns
out to be a detour: here everything can be done 
in terms of line bundles.
More precisely, we have the following
generalization of~\cite[Thm.~1]{BBSw3},
compare also~\cite[Remark, p.~965]{BBSw3}:

\begin{theorem}\label{projcrit}
Let $U \subset X$ be an $N$-invariant open 
subset admitting a good quotient 
$U \to U \quot T$ with $U \quot T$ projective.
Then there is an ample $G$-linearized line bundle 
$L$ on $X$ such that $U = X^{ss}(L,T)$ holds.
Moreover, we have $X = G \mal U$, and $X$ is a 
projective variety.
\end{theorem}

Combining this result 
with~\cite[Thm.~2.1]{Mu} gives
the following supplement to the Hilbert-Mumford 
Theorem~\ref{semisimhimu}:

\begin{corollary}\label{semisimple2hilbmum}
Let $U \subset X$ be as in Theorem~\ref{projcrit}.
Then the intersection $W(U)$ of all translates $g \mal U$, 
$g \in G$, is an open $T$-saturated subset of $U$,
there is a good quotient $W(U) \to W(U) \quot G$, 
and $W(U) \quot G$ is projective. 
\end{corollary}

We come the proof of Theorem~\ref{projcrit}.
A first ingredient is an observation
due to Bia\l ynicki-Birula and \'Swi\c{e}cicka
concerning semisimple group actions on 
the projective space:

\begin{lemma}\label{BBSwlemma}
Let $G$ act on $\PP_{n}$. Then the translates
$g \mal \PP_{n}^{ss}(\mathcal{O}(1), T)$, where $g \in G$,
cover $\PP_{n}$. 
\end{lemma}

\begin{proof} 
Consider the complement $Y$ of the union of all translates
$g \mal \PP_{n}^{ss}(\mathcal{O}(1), T)$, where $g \in G$.
Then $Y$ is empty, because otherwise \cite[Lemma, p.~963]{BBSw3} 
would provide a $T$-semistable point in some irreducible 
component of $Y$.
\end{proof}

The second ingredient of the proof is the following
refinement of Sumihiro's Embedding Theorem, 
compare~\cite[Thm.~1]{Su} and~\cite[Prop.~1.7]{Mu}:

\begin{lemma}\label{embed}
Let $D$ be a $G$-linearized Cartier divisor.
If $X = G \mal X^{ss}(D,T)$ holds,
then there is a $G$-equivariant locally closed
embedding $X \subset \PP_{n}$ such that 
$X^{ss}(D,T)$ is $T$-saturated in
$\b{X} \cap \PP_{n}^{ss}(\mathcal{O}(1), T)$,
where $\b{X}$ is the closure of $X$ in~$\PP_{n}$.
\end{lemma}

\begin{proof}
Let $\mathcal{A}$ be the graded $\mathcal{O}_{X}$-algebra
associated to $D$, and let $U := X^{ss}(D,T)$. 
Since we assumed $X = G \mal U$, Prop.~\ref{ampleinvar} 
tells us that the divisor $D$ is in fact ample.
Moreover, replacing $D$ with a multiple, 
we may even assume that $D$ is very ample, 
and that there are $T$-invariant
$f_{1}, \ldots, f_{r} \in \mathcal{A}_{D}(X)$
such that the sets $X \setminus Z(f_{i})$ are 
affine and cover $U$.

Choose any $G$-invariant vector subspace 
$M \subset \mathcal{A}_{D}(X)$ of finite dimension 
such that $f_{1}, \ldots, f_{r} \in M$ holds, and the 
corresponding morphism $\imath \colon X \to \PP(N)$ is a 
locally closed embedding, where $N$ is the dual $G$-module 
of $M$. 
Then $\imath$ is $G$-equivariant, and $\mathcal{A}_{D}$ equals 
as a $G$-sheaf the pullback of $\mathcal{O}(1)$.
Moreover, by construction, the $f_{i}$ extend to 
$T$-invariant sections of $\mathcal{O}(1)$.
\end{proof}

\begin{proof}[Proof of Theorem~\ref{projcrit}]
First note that $U$ is as well a qp-maximal $N$-set.
Thus we can choose a $G$-linearized Weil divisor $D$ on 
$X$ as in Prop.~\ref{max2semistab}~(i).
By Prop.~\ref{ampleinvar}, $D$ is an ample Cartier
divisor on $X_0 := G \mal U$.
In particular, on $X_{0}$ the $G$-sheaf $\mathcal{A}_{D}$ 
is the sheaf of sections of a $G$-linearized line bundle.

Now choose a locally closed $G$-equivariant embedding 
$X_{0} \subset \PP_{n}$ as in Lemma~\ref{embed}, and
let $\b{X}_{0}$ denote the closure of $X_{0}$ in $\PP_{n}$.
Since $U \quot T$ is complete, we obtain
$$ U = \b{X}_{0} \cap \PP_{n}^{ss}(\mathcal{O}(1), T). $$
Moreover, from Lemma~\ref{BBSwlemma} we infer that the translates
$g \mal U$, where $g \in G$, cover $\b{X}_{0}$.
But this means that we have $X_0 = \b{X}_{0}$.
In particular $X_{0}$ is projective, $X = X_{0}$ holds,
and $D$ is ample.
\end{proof}

\bibliography{}

\end{document}